\documentclass{article}
\usepackage[utf8]{inputenc}
\usepackage[]{amsfonts}
\usepackage{amssymb,latexsym,amsmath,epsfig,mathrsfs}
\usepackage{amsthm}
\usepackage[all]{xy}
\usepackage{graphicx}
\usepackage{float}
\usepackage{textcomp}
\usepackage{fancyhdr}
\usepackage{bm}
\usepackage{mathtools}
\usepackage{tikz}
\usetikzlibrary{matrix,arrows,decorations.pathmorphing}
\usepackage{tikz-cd}
\usepackage{listings}
\usepackage{fullpage}
\usepackage{microtype}
\usepackage{colonequals}
\usepackage{import}
\usepackage{caption}
\usepackage{hyperref}
\hypersetup{colorlinks = false}

\DeclareMathOperator{\im}{Im}

\DeclareMathOperator{\Id}{Id}

\DeclareMathOperator{\GL}{GL}
\DeclareMathOperator{\SL}{SL}

\DeclareMathOperator{\col}{Col}

\DeclareMathOperator{\Ann}{Ann}

\DeclareMathOperator{\rk}{rk}
\DeclareMathOperator{\rel}{rel}
\DeclareMathOperator{\area}{area}

\DeclareMathOperator{\Avg}{Avg}
\DeclareMathOperator{\proj}{proj}

\newcommand{\Q}{\mathbb{Q}}
\newcommand{\R}{\mathbb{R}}
\newcommand{\C}{\mathbb{C}}

\newcommand{\cH}{\mathcal{H}}

\newcommand*\bb[1]{\mathbb{#1}}
\renewcommand*\cal[1]{\mathcal{#1}}

\newcommand*\clos[1]{\overline{#1}}

\newcommand*\define[1]{\textbf{#1}}
\title{Quasidiagonals in Strata of Translation Surfaces}
\author{Christopher Zhang}
\date{\today}

\begin{document}

\theoremstyle{plain}
\newtheorem{theorem}{Theorem}[section]
\newtheorem{thm}[theorem]{Theorem}
\newtheorem{prop}[theorem]{Proposition}
\newtheorem{conj}[theorem]{Conjecture}
\newtheorem{lemma}[theorem]{Lemma}

\newtheorem{prob}[theorem]{Problem}
\newtheorem{axiom}[theorem]{Axiom}
\newtheorem{cor}[theorem]{Corollary}

\newtheorem{alg}[theorem]{Algorithm}
\newtheorem{form}[theorem]{Formula}

\theoremstyle{remark}
\newtheorem{rmk}[theorem]{Remark}
\newtheorem{ques}[theorem]{Question}
\newtheorem{motto}[theorem]{Motto}
\newtheorem{claim}[theorem]{Claim}
\newtheorem{subclaim}[theorem]{Subclaim}
\newtheorem{assumption}[theorem]{Assumption}
\newtheorem{question}[theorem]{Question}

\theoremstyle{definition}
\newtheorem{defn}[theorem]{Definition}
\newtheorem{ex}[theorem]{Example}

\maketitle


\section{Introduction}
For $\cal H$ a component of a stratum of the moduli space of translation surfaces, there is an $\GL^+(2,\R)$-action on $\cal H$. The breakthrough work of Eskin, Mirzakhani, and Mohammadi in \cite{EM} and \cite{EMM} showed that $\GL^+(2,\R)$ orbit closures of translation surfaces are immersed submanifolds that are cut out by linear equations in period coordinates, and furthermore Filip \cite{Filip} showed that they are also subvarieties. For a multi-component stratum $\cal H_1\times\cdots\times\cal H_n$, we have a diagonal $\GL^+(2,\R)$ action.

\begin{defn}
Let $\cal H = \cal H_1\times\cdots\times \cal H_n$ be a product of
connected components of strata of translation surfaces. An
\define{invariant subvariety} of multi-component surfaces
is a closed $\GL^+(2,\R)$-invariant irreducible variety $\cal L\subset{\cal H}$
that is cut out by linear equations with real coefficients in period coordinate
charts. Invariant subvarieties should be assumed to be single-component unless
otherwise specified. The term invariant subvariety includes whole strata.
\end{defn}

\begin{defn} \label{def:pis}
A invariant subvariety of multi-component surfaces $\cal M$ is called a
\define{prime invariant subvariety} if it cannot be written as a product 
$\cal M_1\times \cal M_2$ of two other multi-component invariant subvarieties. 
\end{defn}

One source of these multi-component invariant subvarieties comes from the 
WYSIWYG boundary of Mirzakhani and Wright \cite{MirzakhaniWright_Boundary}.
Starting with a stratum of single-component translation surfaces, going to the
boundary may produce multi-component surfaces. Chen and Wright proved in \cite[Theorem 1.2]{ChenWright_WYSIWYG} that the boundary of invariant subvarieties are multi-component invariant subvarieties. 

\begin{lemma}
Let $\cal M$ be an invariant subvariety. Then the \define{diagonal} in $\cal M\times \cal M$ defined as 
\begin{align*}
    \cal{D}:=\{(M,M)\in \cal M \times \cal M : M\in \cal M\}
\end{align*}
is an invariant subvariety. The \define{antidiagonal} defined as
\begin{align*}
    \overline{\cal{D}}:=\{(M,-\Id(M))\in \cal M\times \cal M:M\in \cal M\}
\end{align*}
is also an invariant subvariety.

\end{lemma}
\begin{proof}
$\cal D$ is $GL^+(2,\R)$-invariant since the action is the diagonal action. Let $a_1,\dots,a_n$ be periods of the first surface while $b_1,\dots,b_n$ be the periods of the same saddles on the second surface. Then, $a_i=b_i$ are equations that cut out $\cal D$. The proof for $\overline{\cal D}$ is similar.
\end{proof}

\begin{rmk}
In hyperelliptic strata, the diagonal and antidiagonal are the same.
\end{rmk}

We generalize the above examples in the following definition.

\begin{defn}
Let $\cal M_1, \cal M_2$ be invariant subvarieties. A prime invariant subvariety $\Delta$ is a \define{quasidiagonal} in $\cal M_1\times \cal M_2$ if the projection maps $\proj_i:\Delta\to \cal M_i$ are dominant. As a shorthand, we will write ``$\Delta \subset \cal M_1\times  \cal M_2$ is a quasidiagonal". By Proposition \ref{ergodic} and Remark \ref{areas} below, every quasidiagonal has a corresponding quasidiagonal where both sides have equal area. Thus, we will assume throughout the paper that both components have equal area unless otherwise stated.
\end{defn}

\begin{rmk}
Any prime invariant subvariety $\Delta$ is a quasidiagonal in 
$\overline{\proj_1(\Delta)}\times \overline{\proj_2(\Delta)}$.
\end{rmk}

\begin{defn}
    Let $(X,\omega)$ is a \define{hyperelliptic translation surface} if there exists a hyperelliptic involution $j:X\to X$ such that $j^*(\omega) = -\omega$. Let $\cal H$ be a component of a stratum of translation surfaces. The \define{hyperelliptic locus} in $\cal H$ is an invariant subvariety that consists of all hyperelliptic translation surfaces of $\cal H$.
\end{defn}

The following is the main theorem of the paper and proves a conjecture by Apisa and Wright \cite[Conjecture 8.35]{ApisaWright_Diamonds} in the case of Abelian differentials.

\begin{thm}\label{main_theorem}
Let $\cal M_i$ be either a connected components of stratum of translation 
surface in genus at least $2$ (without marked points) or a hyperelliptic locus
in such a stratum for $i=1,2$. There exists (equal area) quasidiagonals
$\Delta \subset \cal M_1\times \cal M_2$ only if $\cal M_1 = \cal M_2$.
In this case, $\Delta$ must be the the diagonal or antidiagonal.
\end{thm}

\begin{ex} \label{ex:marked-points}
We do not allow marked points because this will give rise to many uninteresting examples of quasidiagonals. For example, let $\Delta\subset \cal M_1\times\cal M_2$ be a quasidiagonal. Then $\{((M_1,p),M_2):(M_1,M_2)\in \Delta, p\in M_1\}$ is a quasidiagonal in $\cal M_1^{*}\times \cal M_2$. (Here $\cal M_1^*$ denotes the invariant subvariety which is $\cal M_1$ along with a free marked point.)
\end{ex}

Classifying quasidiagonals is helpful with inductive arguments that use the WYSIWYG boundary. In addition, this classification is interesting since quasidiagonals show relationships between $\cal M_1$ and $\cal M_2$. 

\begin{defn}
A continuous, $SL(2,\R)$-invariant map between invariant subvarieties $\phi: \cal M\to \cal N$ is called a \define{morphism} if it is linear in period coordinates. 
\end{defn}

A morphism $\phi:\cal M\to \cal N$ between invariant subvarieties gives a quasidiagonal $\{(M,\phi(M)): M\in \cal M\}\subset \cal M\times \overline{\phi(\cal M)}$. For example, when $\cal H$ is not hyperelliptic, $-\Id$ gives rise to a nontrivial automorphism of $\cal H$. This corresponds to the antidiagonal.  

\begin{cor}
Let $\cal H,\cal H'$ be strata. There are no dominant morphisms $\phi: \cal H\to \cal H'$ other than $\Id,-\Id:\cal H\to \cal H$. 
\end{cor}

The above support the following heuristic: a quasidiagonal $\Delta\subset \cal M_1\times \cal M_2$ exists if and only if $\cal M_1$ and $\cal M_2$ are ``related''. 

\begin{ex} \label{ex:branched-covers}
    Let $\widetilde {\cal H}(2,0^2) \subset \cal H(2^2,1^2)$ be the space of all double covers of surfaces in $\cal H(2)$ branched at two marked points. There is a quasidiagonal $\Delta\subset \cal H(2)\times \widetilde{\cal H}(2,0^2)$ consisting of all $(M,\widetilde M)$ where $\widetilde M$ is a branched double cover of $M$.
\end{ex}

\begin{lemma}
    Given quasidiagonals $ \Delta_L \subset \cal M_1\times \cal M_2$ and $\Delta_R \subset \cal M_2\times \cal M_3$, then $\Delta_L * \Delta_R \subset \cal M_1\times \cal M_3$ defined as the closure of $\{(M_1,M_3)\in \cal M_1\times \cal M_3:\exists M_2 \text{ such that } (M_1,M_2)\in \Delta_L,(M_2,M_3)\in \Delta_R\}$ is a quasidiagonal.
\end{lemma}
\begin{proof}
    By Theorem \ref{ChenWright}, for $(M_1,M_2)\in \Delta_L$
    the absolute periods of $M_1$ locally determine the absolute periods of $M_2$, and
    similarly for $(M_2,M_3)\in \Delta_R$, the absolute periods of $M_2$ locally
    determine the absolute periods of $M_3$. Thus, the absolute periods of 
    each side of $\Delta_L*\Delta_R$ locally determines the absolute periods
    of the other. 
    Thus, $\Delta_L*\Delta_R$ cannot be the product of invariant subvarieties
    of multi-component surfaces, so it must be prime.
    $\proj_2(\Delta_L)\cap \proj_1(\Delta_R)$ is a Zariski open subset of
    $\cal M_2$. $\proj_1: \Delta_L*\Delta_R\to \cal M_1$ is dominant since
    $\proj_1(\Delta_L * \Delta_R) \supset \proj_1(\proj_2^{-1}(\proj_2(\Delta_L)
    \cap \proj_1(\Delta_R)))$.
\end{proof}

\begin{cor}
    Let $\cal M_1,\cal M_2$ be invariant subvarieties. We define $\cal M_1\sim \cal M_2$ if there exists a quasidiagonal $\Delta \in \cal M_1\times \cal M_2$. Then, $\sim $ is an equivalent relation.
\end{cor}

Theorem \ref{main_theorem} implies that there are many distinct $\sim$ equivalence classes. It would be interesting to classify these equivalence classes. The only ways the author is aware of to make related invariant subvarieties is through adding marked points as in Example \ref{ex:marked-points} or branched covering constructions as in Example \ref{ex:branched-covers}. 
\begin{question}
    Let $\cal R$ be an equivalence class of related invariant subvarieties with marked points in rank $\ge 2$. Does $\cal R$ contain a minimal element through which all other related invariant subvarieties can be produced through adding marked points and branched covering constructions?
\end{question}

Another application of this work is to measurable joinings of Masur-Veech measures.

\begin{defn}
Let $(X_1, \mu_1, T_1)$ and $(X_2, \mu_2, T_2)$, where $\mu_i$ is a measure on a space $X_i$ and $T_i:X_i\to X_i$ is a measure preserving transformation. A \define{joining} is a measure on $X_1\times X_2$ invariant under the product transformation $T_1\times T_2$, whose marginals on $X_i$ are $\mu_i$. A measure $\mu$ on a space $X$ is \define{prime} if it cannot be written as a product $\mu = \mu_1\times\mu_2$, $X = X_1\times X_2$, where $\mu_i$ is a measure on $X_i$.
\end{defn}

Assuming a multi-component version of Eskin-Mirzakhani's measure classification result, Theorem \ref{main_theorem} classifies ergodic measurable joinings of Masur-Veech measures on strata.

\begin{assumption} \label{magic-wand} (See \cite[Conjecture 2.10]{MirzakhaniWright_Boundary}) We define an affine measure as in \cite[Definition 1.1]{EM}.
For a multi-component stratum under the diagonal action of $\SL(2,\R)$, the only ergodic invariant measures are $\SL(2,\R)$-invariant and affine.
\end{assumption}

\begin{cor}
Let $\cal M_i$ be either a stratum of translation surfaces or the hyperelliptic locus in such a stratum and $\mu_i$ be the Masur-Veech measure on the unit area locus of $\cal M_i$, for $i=1,2$. Under Assumption \ref{magic-wand}, the only prime ergodic joinings of $\mu_1,\mu_2$ (under the diagonal action of $\SL(2,\R)$) are the Masur-Veech measure on the diagonal or antidiagonal.
\end{cor}

\begin{proof}
    By Assumption \ref{magic-wand}, the only ergodic $\SL(2,\R)$-invariant measures are affine, so they are supported on an invariant subvariety $\cal M$. Since the joining is prime, $\cal M$ must be prime. The condition on the marginals implies that $\cal M$ must be a quasidiagonal in $\cal M_1\times \cal M_2$. Our assumptions also guarantee equal area on both sides. We may now apply Theorem \ref{main_theorem}.
\end{proof}

\textbf {Organization:} Section 2 is the bulk of the paper and contains a summary of the techniques used in the proof of the main theorem. The proof will use induction. Section 3, which is slightly technical, is the base case of the proof. The heart of the proof is in Section 4. Section 5 is an appendix.
\vspace{3mm}

\textbf {Acknowledgements:} During the preparation of this paper, the author was partially supported by NSF Graduate Research Fellowship DGE 1841052. The author is grateful to his advisor Alex Wright for suggesting the problem and for significant help and guidance. He is also grateful to Paul Apisa for helpful discussions pertaining to this problem. 

\section{Background}
In this section, we define notation and terminology and list the background needed in the paper. We also prove the results needed in the proof of the main theorem.
By ``stratum", we will refer to a ``connected component of a stratum of connected translation surfaces without marked points in genus at least $2$" unless otherwise stated. Let $\cal H$ be a stratum. We will often use a single letter $M = (X,\omega)\in \cal H$ to denote a translation surface, where $X$ is the underlying Riemann surface and $\omega$ the holomorphic $1$-form. 

\subsection{Prime Invariant Subvarieties}

The notion of a prime invariant subvariety (Definition \ref{def:pis}) was defined and studied in \cite{ChenWright_WYSIWYG}. Here we list some results about them.

\begin{thm}[Chen-Wright] \label{ChenWright} \cite[Theorem 1.3]{ChenWright_WYSIWYG} 
    For a prime invariant subvariety 
    $\cal P \subset \cal H_1\times \cdots\times \cal H_n$, the absolute periods
    locally determine each other. 
\end{thm}

\begin{defn} \label{def:rank}
    Let $\cal M$ be a single-component invariant subvariety with marked points,
    and $M\in \cal M$. 
    Let $p:H^1(M,\Sigma;\C)\to H^1(M;\C)$ be the forgetful map. The \define{rank} of $\cal M$ is $\frac{1}{2}\dim p(T_M\cal M)$.
    The \define{rank} of a prime invariant subvariety $\cal N$ is 
    $\rk (\overline{\proj_i(\cal N)})$, which is independent $M$,
    and it is independent of $i$ by 
    Theorem \ref{ChenWright}. $\cal M$ is \define{full rank} if $\rk \cal M$
    is equal to the rank of $\cal H$, the stratum 
    that contains $\cal M$. 
\end{defn}

\begin{rmk} 
    By \cite[Theorem 1.4]{AEM}, $p(T_M \cal M)$ is symplectic, so the rank of an invariant subvariety is always an integer.
\end{rmk}

\begin{prop}[Chen-Wright] \label{ergodic} \cite[Corollary 7.4]{ChenWright_WYSIWYG}
In a prime invariant subvariety, the $g_t$ action is ergodic on the unit area locus. Thus, the ratio of areas of the components is constant. 
\end{prop}

\begin{rmk} \label{areas}
For any quasidiagonal $\Delta$ we get an infinite number of quasidiagonals $\Delta_r = \{(M_1,rM_2):(M_1,M_2)\in \Delta\}$. 
By Proposition \ref{ergodic}, we can scale $\Delta$ so that both components have the same area.
\end{rmk}

\subsection{Cylinder Deformations}

By a \define{cylinder} $C$ we refer to a maximal topological open annulus foliated by closed geodesics. If $C$ is a cylinder on $M$, there is a corresponding cylinder on every surface in a small enough neighborhood of $M$. As an abuse of notation, we often refer to these corresponding cylinders as $C$. By \define{core curve} of a cylinder we refer to one of these closed geodesics. A cylinder has two \define{boundary components}, which may intersect at singularities. A \define{cross curve} is a saddle connection that goes from one boundary component to the other.

Let $\cal M \subset \cal H = \cal H_1\times\dots \times \cal H_n$ be
an invariety subvariety of multi-component surfaces.
We reproduce a theorem by 
Smillie-Weiss \cite[Theorem 5]{SmillieWeiss}. 
The original theorem was only proven in the single-component case. The proof
of the multicomponent case has been deferred to the appendix.
Let $H_t = \left\{
\begin{pmatrix}
1 & t \\
0 & 1
\end{pmatrix}
: t\in \R
\right\} \subset \SL(2,\R)$.

\begin{thm}[Smillie-Weiss 2004] \label{SmillieWeiss}
Every $H_t$ orbit closure in a stratum of multicomponent quadratic 
differentials contains a surface $q$ such that every component of $q$ is 
horizontally periodic.
\end{thm}

\begin{lemma}\label{real}
Let $M=(M_1,\dots,M_n)\in \Delta$ be a surface in a prime invariant subvariety. 
Choose a period coordinate chart $U$ around $M$ and let 
$M^t=(M_1^t,\dots,M_n^t)$, $t\in[0,1]$, be a path in $U$ with $M^0=M$. If for some $i$, for each $t$ the imaginary parts of the absolute periods of $M^t_i$ are the same, then the same is true for all $i$. 
\end{lemma}
\begin{proof}
By Theorem \ref{ChenWright}, the absolute periods of $M_i$ determine each other and $\cal M$ is cut out by equations with real coefficients in period coordinate charts.
\end{proof}

\begin{lemma} \label{horizontally_periodic}
Let $M=(M_1,\dots,M_n) \in \cal M $ be a prime invariant subvariety. If $M_i$ is horizontally periodic, then $M_j$ must be horizontally periodic for every $1\le j \le n$. In this case, we say that the pair $M$ is \define{horizontally periodic}.
\end{lemma}
\begin{proof}
Assume by contradiction and with loss of generality that $M_1$ is horizontally periodic and $M_2$ is not. By Theorem \ref{SmillieWeiss}, there is a sequence $t_n\to \infty$ such that 
\[
\begin{pmatrix}
1 & t_n\\
0 & 1
\end{pmatrix}M_2\to M_2^\infty 
\]
where $M_2^\infty$ is horizontally periodic. Since $M_1$ is horizontally 
periodic, its orbit closure under the action of $H_t$ is an $n$-dimensional torus $T$, which 
contains surfaces with locally the same imaginary periods. Since the $H_t$-orbit closure of $M$ is compact, after passing to a subsequence, we also have that $\begin{pmatrix}
1 & t_n\\
0 & 1
\end{pmatrix}M_1\to M_1^\infty$ and $M_1^\infty$ must also be horizontally periodic. We choose a coordinate chart $U\times V$ centered around $M^\infty = (M_1^\infty,\dots, M_n^\infty)$, where $U\subset \cal M_1$ and $V\subset \cal M_2\times \cdots \times \cal M_n$. For large enough $t_n$, $M^{t_n}_1 := 
\begin{pmatrix}
1 & t_n\\
0 & 1
\end{pmatrix}M_1 \in U$ and since $M^{t_n}_1$ and $M^\infty_1$ all lie on $T$ they all have the same real periods. However, for large enough $n$, there is a surface $M_2^{t_n} := \begin{pmatrix}
1 & t_n\\
0 & 1
\end{pmatrix}M_2$ such that there is some cylinder $C$ of $M^\infty_2$ that persists on $M^{t_n}_2$ but it is not horizontal. Thus, the core curve of $C$ is an absolute period that changes in imaginary part on a path from $M^{t_n}_2$ to $M_2^\infty$. This contradicts Lemma \ref{real}. 
\end{proof}

Let $\cal C = \{C_1,\dots,C_r\}$ be a collection of cylinders on $M = (M_1,\dots,M_n)\in \cal M$ and let $\gamma_i$ be a core curve of $C_i$. The tangent space $T_M\cal M$ is a subspace of $T_M\cal H = H^1(M,\Sigma;\C) := H^1(M_1,\Sigma;\C)\times\cdots \times H^1(M_n,\Sigma;\C)$. Thus, there is a projection $\pi:H_1(M,\Sigma;\C)\to (T_M\cal M)^*$. We view $\gamma_i$ as elements of $H_1(M,\Sigma;\C)$ by setting the $H_1(M_j,\Sigma;\C)$ to be zero on the components $M_j$ that do not contain $C_i$.

\begin{defn} \label{M-parallel}
    $\cal C$ is called \define{$\cal M$-parallel} if all $\pi(\gamma_i)$ are colinear in $(T_M\cal M)^*$. Being $\cal M$-parallel is an equivalence relation on cylinders, so we call $\cal C$ an \define{$\cal M$-parallel class} if it is an equivalence class of $\cal M$-parallel cylinders.
    \end{defn} 
    
    Intuitively, $\cal M$-parallel means that there is a neighborhood
    $M\in U \subset \cal M$ such that all cylinders in $\cal C$ remain parallel in
    this neighborhood. See \cite{Wright_Cylinder} for a more detailed discussion on
    $\cal M$-parallel cylinders.

Now we continue to the statement of the Cylinder Deformation Theorem.
Let $\cal C$ consist of cylinders $C_1,\dots,C_r$, with heights $h_1,\dots,h_r$, and let $\alpha_i$ be the cohomology class associated to the core curve of $C_i$ under Poincare duality $H_1(M-\Sigma;\C)\cong H^1(M,\Sigma;\C)$. Define 
\[
\eta_{\cal C} = \sum_{i=1}^r h_i\alpha_i
\]
See also \cite[Section 2]{Wright_Cylinder} for a more detailed definition of $\eta_{\cal C}$. For horizontal cylinders, moving in the direction of $i\eta_{\cal C}$ in period coordinates stretches all the cylinders in $\cal C$ (in proportion to the height of the cylinder) and doesn't change the rest of the surface. Moving in the direction of $\eta_{\cal C}$ shears all the cylinders of $\cal C$ and doesn't change the rest of the surface. We are now able to state the theorem.

\begin{thm}[Cylinder Deformation Theorem, Wright 2015]\label{cylinder_deformation}
Let $\cal M \subset \cal H$ be an invariant subvariety of a multi-component stratum. Let $\cal C$ be a $\cal M$-parallel class of cylinders on $M$. Then $\eta_{\cal C} \in T_{M}\cal M$. (Here $T_M\cal M$ denotes the tangent space to $\cal M$ at $M$, which is a subspace of $H^1(M, \Sigma;\C)$).  
\end{thm}

We will sketch a proof of Wright's Cylinder Deformation Theorem
for multi-component surfaces because the original theorem
was only stated for single-component surfaces. The proof is identical to
the original proof but we omit many details, which can be found in the original
paper \cite{Wright_Cylinder}.
The following lemma is \cite[Lemma 3.1]{Wright_Cylinder}.
\begin{lemma} \label{linear-flows-torus}
Let $M$ be horizontally periodic and $\cal C$ be an $\cal M$-parallel class of cylinders. Let the moduli of the cylinders of $\cal C$ be independent over $\Q$ of the moduli of the remaining horizontal cylinders. Then $\eta_{\cal C}\in T_M\cal M$.
\end{lemma}
\begin{proof}
The $H_t$-flow of a horizontally periodic surface is the same as the flow on a $r$-dimensional torus whose slope is determined by the moduli of the cylinders. 
\end{proof}

The following lemma can be found in \cite[Lemma 4.9]{Wright_Cylinder}.

\begin{lemma}\label{linear-algebra-colinear}
Let $V$ be a finite-dimensional vector space and $F\subset V^*$ a finite collections of linear functionals on $V$, no two of which are colinear. The collection of functions $1/w$ for $w\in F$ are linearly independent over $\R$. This remains true when the functions are restricted to any nonempty open set of $V$. 
\end{lemma}

\begin{proof}[Proof of Theorem \ref{cylinder_deformation}]
Let $\cal C$ be an $\cal M$-parallel class of cylinders on $M$. By Theorem 
\ref{SmillieWeiss}, there is a horizontally periodic surface $M'$ in the 
$H_t$ orbit closure of $M$. The corresponding set of cylinders, which we still call $\cal C$, is an $\cal M$-parallel class on $M'$. 

\begin{claim}
There is a surface $M''$ that is a real deformation of $M'$ such that the moduli of the cylinders in $\cal C$ are independent over $\Q$ of the moduli of the cylinders not in $\cal C$. \cite[Lemma 4.10]{Wright_Cylinder}
\end{claim}

By the claim we have a surface $M''$ where the moduli of the cylinders in $\cal C$ are independent of the rest of the cylinders. Lemma \ref{linear-flows-torus} finishes the proof. Thus, it suffices to prove the claim. Let $C_{r+1},\dots,C_l$ be the cylinders of $M'$ not in $\cal C$, and let $m_i$ be the moduli of the cylinder $C_i$. Assume by contradiction that the claim is false. Because $\cal M$ is cut out by real linear equations in period coordinates, there is some rational relation that holds for small real deformations of $M'$
\[
\sum_{i=1}^r q_im_i = \sum_{j=r+1}^{l} q_jm_j
\]
for some $q_i\in \Q$, where neither the right hand nor the left hand side is identically zero. Recall that $m_i = h_i/c_i$ is the modulus of a cylinder. Since the cylinders in $\cal C$ are all $\cal M$-parallel, the $c_i$ are all multiples of each other in a small neighborhood. Allowing the coefficients to be real numbers, we can remove all but one representative from each $\cal M$-parallel class. Thus, we have an equation of the form
\[
r_1m_1 = \sum_{\substack{j\in J, \\ J\subset \{r+1,\dots,l\}}} q_jm_j
\]
where neither side is zero, no two of the cylinders $C_{j}$ are $\cal M$-parallel. However, $1/m_{i_j}$ are linear functional (over an open set in the space of real deformations) that are not colinear, so this relation cannot hold by Lemma \ref{linear-algebra-colinear}. This is a contradiction, so the claim and the theorem are proven. 
\end{proof}

The following corollary, which is \cite[Proposition 3.2]{NW}, immediately generalizes to the multi-component setting:

\begin{cor}\label{cylinder_proportion}
    Let $M$ be a surface in an invariety subvariety of multi-component surfaces
    $\cal M$. Let $\cal C,\cal C'$ be $\cal M$-parallel classes of cylinders on $M$, and let $C,D$ be cylinders in $\cal C'$. Then,
    \[
        \frac{\area(C\cap \cal C)}{\area(C)} = \frac{\area(D\cap \cal C)}{\area(D)}.
    \]
\end{cor}

 

\begin{lemma} \label{m-parallel-quasidiagonal}
    Let $\cal M$ be a prime invariant subvariety and $\cal C$ be a
    $\cal M$-parallel class of cylinders on a surface 
    $M=(M_1,\dots,M_n)\in \cal M$. Define $\cal M_i := \overline{\proj_i(\cal M)}$
    to be the closure of the $i$-th projection of $\cal M$. Let $\cal C_i$ be
    the cylinders of $\cal C$ on $M_i$. Then, $\cal C_i$ is a nonempty
    $\cal M_i$-parallel class of cylinders. Furthermore, if $\cal C_i, \cal C_i'$ are two distinct $\cal M_i$-parallel classes, then there are distinct $\Delta$-parallel classes $\cal C,\cal C'$ that contain $\cal C_i,\cal C_i'$ respectively.
\end{lemma}

\begin{proof}
    First we show each $\cal C_i$ is nonempty. Assume by contradiction that $M_i$ does not have a cylinder in $\cal C$, but $M_j$ does. By the Cylinder Deformation Theorem, we can perform standard cylinder dilation on $\cal C$ while remaining in $\Delta$. This causes the absolute periods of $M_j$ to change without changing the absolute periods of $M_i$, which contradicts Theorem \ref{ChenWright}.

    Let $p:H^1(M,\Sigma;\C) \to H^1(M;\C) $ be the projection from relative to absolute cohomology. Then, $(pT_M\cal M)^*\subset (T_M\cal M)^*$. Since the core curves of cylinders $\gamma_i$ are elements of absolute homology $H_1(M;\C)$, we have $\pi(\gamma_i) \in (pT_M\cal M)^*$ (where $\pi$ is defined in the discussion before Definition \ref{M-parallel}).
    Now, let $\{\gamma_j\}$ be the core curves of cylinders of $M_i$. By \cite[Theorem 1.3]{ChenWright_WYSIWYG}, $(pT_M\Delta)^* \cong (pT_{M_i}\cal M_i)^*$, so $\pi(\gamma_i)$ are colinear in $(pT_{M_i}\cal M_i)^*$ if and only if they are colinear in $(pT_M\Delta)^*$. Thus, the $\gamma_i$ they are $\cal M_i$-parallel if and only if they are $\Delta$-parallel. This proves the rest of the lemma.
\end{proof}



\subsection{WYSIWYG Compactification}

We give a short overview of the WYSIWYG compactification. See \cite{MirzakhaniWright_Boundary} and \cite{ChenWright_WYSIWYG} for more formal introductions.

\begin{defn}
Let $\cal H,\cal H'$ be strata of multi-component translation surfaces potentially having marked points. Let $M_n = (X_n,\omega_n)\in \cal H$ and $\Sigma_n$ be its set of singularities and marked points, and let $M = (X,\omega) \in \cal H'$ and $\Sigma$ its set of singularities and marked points. We say that $M_n$ converges to $M$ if there are decreasing neighborhoods $\Sigma \subset U_i \subset M$ such that there are $g_i: X-U_i\to X_i$ that are diffeomorphisms onto their images satisfying
\begin{enumerate}
    \item $g_i^*(\omega_i) \to \omega$ in the compact-open topology on $M-\Sigma$.
    \item The injectivity radius of points not in the image of $g_i$ goes to zero uniformly in $i$.
\end{enumerate}
See \cite[Definition 2.2]{MirzakhaniWright_Boundary}.
\end{defn}

Thus, we can construct $\partial \cal H$ from $\cal H$ by including all $\cal H'$ such that a sequence of surfaces in $\cal H$ converges to a surface in $\cal H'$. Multiple copies of a stratum can be included if there are two sequences that converge to the same surface in $\cal H'$ but are not close in $\cal H$. We call the union $\clos{\cal H} = \cal H \cup \partial \cal H$ (with the topology given by the above convergence of sequences) the WYSIWYG partial compactification of $\cal H$. For any invariant subvariety $\cal M\subset \cal H$, we define $\partial M$ to be $\overline{\cal M}-\cal M$, where the closure $\overline {\cal M}$ is taken in $\overline {\cal H}$.


\begin{rmk}
Even if $M_n$ is a convergent sequence of surfaces without marked points, its limit may have marked points. 
\end{rmk}

Let $M_n = (X_n,\omega_n)\in \cal M$ be a sequence of multi-component translation surfaces that has a limit $M = (X,\omega) \in \partial \cal M$. Let $\cal H'$ be the stratum with marked points that contains $M$. Let $\cal N$ be the connected component of $\cal H'\cap \partial \cal M$ that contains $M$. We call $\cal N$ the \define{component of the boundary of $\cal M$} that contains $M$. The sequence $X_n$ will approach a limit $X'$ in the Deligne-Mumford compactification. For large enough $n$ there is a map $f_n:X_n\to X'$ called the collapse map. There is also a map $g:X\to X'$ identifying together marked points of $X$. Define $(f_n)_*:H_1(X_n,\Sigma_n)\to H_1(X,\Sigma)$ and $V_n = \ker ((f_n)_*)$. 

\begin{prop}\label{tangent-space-boundary}
After identifying $H_1(X_n,\Sigma_n)$ for different $n$, $V_n$ eventually becomes constant which we call $V$. For large enough $n$, $T_M\cal H'$ can be identified with $\Ann(V)$.
\end{prop}
This proposition was proven for multi-component surfaces in
\cite[Proposition 2.5 and Proposition 2.6]{MirzakhaniWright_Boundary}.


\begin{thm}[Mirzakhani-Wright 2017, Chen-Wright 2021] \label{thm:boundary}
Let $\cal M$ be an invariant variety in a stratum $\cal H$ of connected 
translation surfaces with marked points.
Let $M_n\in \cal M$ be a sequence that converges to $M \in \partial M$.
Let $\cal H'$ be the stratum that contains $M$ and $\cal M'$ be the component
of the boundary of $\cal M$ that contains $M$.
By Proposition \ref{tangent-space-boundary} we identify $T_M\cal H'$
with $\Ann(V)$. 
Then, $T_M \cal M'$ can be identified with $T_{M_n}\cal M \cap \Ann(V)$.
\end{thm} 

This theorem was proven in \cite[Theorem 1.1]{MirzakhaniWright_Boundary} when $M$ is connected, and in \cite[Theorem 1.2]{ChenWright_WYSIWYG} when $M$ is disconnected.



\begin{thm} \cite[Theorem 1.1]{MirzakhaniWright_Full} \label{thm:full_rank}
    An invariant subvariety is full rank if and only if it is a full stratum or
    a hyperelliptic locus. 
\end{thm}

We clarify our definition of hyperelliptic locus in the case of marked points.

\begin{defn} \label{def:hyperelliptic}
    Let $\cal H$ be a stratum of genus $\ge 2$
    with marked points.
    $M\in \cal H$ is called \define{hyperelliptic} if there is an involution
    $J:M\to M$ that is the hyperelliptic involution on the underlying
    Riemann surface of $M$ such that $J$ maps marked points to marked points.
    $\cal H$ is \define{hyperelliptic} if every surface in $\cal H$ is 
    hyperelliptic.
    An invariant subvariety $\cal M\subset \cal H$ is
    a \define{hyperelliptic locus} if it contains exactly the 
    hyperelliptic surfaces 
    of $\cal H$.
\end{defn}

\begin{cor} \label{collapse-full-rank}
    If $\cal M$ is a full rank invariant subvariety, then any component 
    $\cal N$ of
    $\partial \cal M$ is a full rank invariant subvariety with marked points.
\end{cor}
\begin{proof}
    By Theorem \ref{thm:full_rank}, a full rank invariant subvariety is a
    stratum or hyperelliptic locus.
    A stratum is cut out by no equations in period coordinates, so
    by Theorem \ref{thm:boundary}
    $\cal N$ is cut out by no equations, so it
    is a stratum.
    A hyperelliptic locus with hyperelliptic involution $J$
    is cut out by the equations
    $\int_\gamma \omega + \int _{J(\gamma)}\omega = 0$
    for all saddles $\gamma$.
    $J$ restricts to the hyperelliptic involution on $\cal N$ and the 
    equations that cut out $\cal M$ restrict to the corresponding equations
    that cut out $\cal N$. Thus, $\cal N$ is a hyperelliptic locus.
\end{proof}

\subsection{Cylinder Collapse} \label{sec:collapse}

We will define cylinders collapse and diamonds in a similar fashion to \cite[Lemma 4.9]{ApisaWright_HighRank}.
Let $ \cal M$ be an invariant subvariety of multi-component surfaces with marked points. Let $M\in \Delta$ and $\cal C$ an $\cal M$-parallel class of horizontal cylinders on $M$. Fix a cross curve $\gamma$ of some cylinder in $\cal C$. We now define the operations $u_s^{\cal C}$ and $a_t^{\cal C}$, which are shearing and scaling the cylinders $\cal C$ respectively.
Let $u_s^{\cal C}(M)$ be the surface obtained by adding $s\eta_{\cal C}$ to $M$ in period coordinates, where $\eta_{\cal C}$ is defined in the paragraph above Theorem \ref{cylinder_deformation}. Define $a_t^{\cal C}(M)$ to be the surface obtained by adding $(e^{t}-1)\eta_{\cal C}$ to $M$ in period coordinates. We define $\col _{\cal C,\gamma} M$ to be the following operation: Choose $s$ so that $\gamma$ is vertical on $u_s^{\cal C}(M)$. Then
\[
    \col _{\cal C,\gamma} M : = \lim_{t\to -\infty} a_t^{\cal C}u_s^{\cal C}M.
\]
This limit exists in the WYSIWYG compactification 
see \cite[Lemma 3.1]{MirzakhaniWright_Boundary} or 
\cite[Lemma 4.9]{ApisaWright_HighRank}. 
This operation has been defined so that $\gamma$ degenerates. If $M$ is a multicomponent surface, this means that the component that contains $\gamma$ degenerates, while the other components of $M$ may or may not degenerate.
We also define
\[
\col _{\cal C,\gamma}\cal M
\]
to be the component of $\partial \cal M$ that contains
$\col _{\cal C,\gamma} M $.
We define $\col_{\cal C,\gamma}\cal C$
to be the $\cal M$-parallel collection of saddles on $\col_{\cal C,\gamma} M$
from the collapsed cylinders in $\cal C$.

Let $\cal C_1,\cal C_2$ be disjoint equivalence classes of cylinders on $M$,
and
let $\gamma_1,\gamma_2$ be cross curves of cylinders of $\cal C_1,\cal C_2$
respectively.
As an abuse of notation, for $j\neq i$, we also use the notation $\cal C_j,\gamma_j$
to refer to the equivalence class of cylinders and saddle connection
on $\col_{\cal C_i,\gamma_i}M$ that correspond to $\cal C_j,\gamma_j$
respectively.
We define
\[
\col_{\cal C_1,\cal C_2,\gamma_1,\gamma_2}M 
= \col_{\cal C_1,\gamma_1}\col _{\cal C_2,\gamma_2}M 
= \col_{\cal C_2,\gamma_2}\col_{\cal C_1,\gamma_1}M.
\]

\begin{defn}
    Let $\cal M$ be an invariety subvariety of multi-component surfaces. A surface $M\in \cal M$ is called $\cal M$-\define{generic} if two saddles on the same component of $M$ are parallel only if they are $\cal M$-parallel. If $\cal M$ is clear from context, we will just call $M$ generic.
\end{defn}

\begin{lemma} \label{generic}
    For any multi-component invariant subvarity, a dense $G_\delta$ set of surfaces are generic.
\end{lemma}
\begin{proof}
The condition that two saddles that are not generically parallel are parallel defines a linear subspace in period coordinates. There are countably many saddles on a surface, and the coefficients of the equation must be in a the field of definition of the component that the two saddles are on. The field of definition is a finite extension of $\bb Q$ by \cite[Theorem 1.1]{Wright_Field}. 
\end{proof}

\begin{lemma} \label{simple-cylinder}
    Let $\cal M$ be a full rank invariant subvariety,
    and let $M\in \cal M$ be generic. Then, every cylinder $C$ on $\cal M$ is simple.
\end{lemma}
\begin{proof}
    By Theorem \ref{thm:full_rank}, $\cal M$ is a stratum or hyperelliptic 
    locus.
    Assume by contradiction there was a cylinder $C$ that was not simple.
    Then, there would be a boundary component with more than one saddle.
    These saddles would be $\cal M$-parallel because $M$ is generic.
    That is not possible in a stratum,
    so $\cal M$ must be a hyperelliptic locus.
    Let $\gamma_1,\gamma_2$ be parallel saddle on the boundary of $C$.
    Every cylinder on $M$ is either fixed or swapped with another cylinder.
    In either case, the quotient surface $N$ has a corresponding cylinder,
    which we still call $C$, with two saddles,
    which we still call $\gamma_1$ and $\gamma_2$.
    Now we appeal to \cite{MasurZorich},
    and we will use the formulation and terminology of
    \cite[Proposition 4.4]{ApisaWright_MarkedPoints}.
    By this theorem, removing $\gamma_1,\gamma_2$ from $N$ disconnects the
    surface into two components $\cal A, \cal B$,
    and the component $\cal B$ not containing $C$ has trivial holonomy.
    Since $\cal M$ is hyperelliptic, $N$ is genus $0$, so $\cal B$ is
    topologically a cylinder. However, gluing this cylinder back along $\gamma_1,\gamma_2$ would create genus. This is a contradiction since $N$ is genus $0$. Thus, the cylinder $C \subset M$ must have been simple.
\end{proof}


\begin{lemma} \label{collapse-one-dim}
    Let $\Delta\subset \cal M_1\times\cal M_2$ be a quasidiagonal, $M\in \Delta$ a generic surface, $\cal C\subset M$ a $\Delta$-parallel class of cylinders on $M$, and $\gamma$ a cross curve of a cylinder in $\cal C$.
    Let $\Delta := \col_{\cal C,\gamma}\Delta$.
    Then $\dim \Delta' = \dim \Delta - 1$ and 
    $\dim \cal M_i - 1 \le \overline{\proj_i (\Delta')} \le \dim \cal M_i$ for $i=1,2$.
\end{lemma}
\begin{proof}
    Let $M' := \col_{\cal C,\gamma}M$.
    Since $M$ is generic, the space of vanishing cycles is one-dimensional, so by Proposition \ref{tangent-space-boundary}, $\dim \Delta' = \dim \Delta - 1$, and $T_{M'} \Delta'$ can be viewed as a one-dimensional subspace of $T_M\Delta$. Then, $\dim (\proj_i T_M\Delta) -1 \le \dim(\proj_i T_{M'}\Delta') \le \dim (\proj_i T_M \Delta)$. The result follows.
\end{proof}

\begin{lemma}\label{lem:codimension_one_IS}
    Let $\cal M$ be an invariant subvariety and $\cal N \subset \cal M$ be a codimension $1$ subvariety. Then $\cal M,\cal N$ must have the same rank.
\end{lemma}
\begin{proof}
    $p(T \cal N)\subset p(T\cal M)$ is codimension at most $1$ and the symplectic form on $p(T\cal M)$ restricts to a symplectic form on $p(T\cal N)$
    , so in fact $p(T \cal M)\cong p(T\cal N)$.
\end{proof}

\begin{lemma} \label{drop_rank_lemma}
Let $\Delta\subset \cal M_1\times\cal M_2$ be a quasidiagonal,
where $\cal M_i$ is a full rank invariant subvariety, $M = (M_1,M_2)\in \Delta$ is generic, $\cal C$ a $\Delta$-equivalence class of cylinders on $M$, and $\gamma$ a cross curve of a cylinder $C\in \cal C$ on $M_1$. If $\col _{C, \gamma} \cal M_1$ is lower rank than $\cal M_1$, then there must be a $\gamma' \subset M_2$ generically parallel to $\gamma$ that is a cross curve of a cylinder in $\cal C$.
\end{lemma}
\begin{proof}
Assume by contradiction no saddle of $\cal M_2$ collapses in
$\col_{\cal C,\gamma}\Delta$, so
$\cal M_2' := \overline{\proj_2(\col_{\cal C, \gamma}\Delta)}\subset \cal M_2$.
However, $\cal M_2'$ is dimension at most one less than $\cal M_2$ by 
Lemma \ref{collapse-one-dim} and has lower rank than $\cal M_2$ by
Theorem \ref{ChenWright}.
This contradicts Lemma \ref{lem:codimension_one_IS}.
Thus, some $\gamma' \subset M_2$ collapses in $\col _{\cal C,\gamma}\Delta$.
A priori $\gamma'$ may cross multiple adjacent cylinders.
Since we assumed $M$ is generic,
by Lemma \ref{simple-cylinder} all cylinders must be simple.
Adjacent simples cylinders meet at marked points,
but we assumed that there are no marked points,
so there are no adjacent cylinders.
\end{proof}

\begin{lemma} \label{full-rank-sides}
    Let $\Delta\subset \cal M_1\times\cal M_2$ be a quasidiagonal,
    where $\cal M_i$ is a full rank invariant subavariety,
    $M = (M_1,M_2)\in \Delta$ is generic,
    $\cal C$ a $\Delta$-equivalence class of cylinders on $M$,
    and $\gamma$ a cross curve a cylinder $C\in \cal C$.
    Let $\cal M_i' := \clos{\proj_i(\col_{\cal C, \gamma} \Delta)}$.
    Then, $\col_{\cal C,\gamma}\Delta \subset \cal M_1' \times \cal M_2'$ is a
    quasidiagonal, and $\cal M_i'$ is a full rank invariant subvariety
    with marked points.
\end{lemma}
\begin{proof}
    By \cite[Lemma 9.1]{ApisaWright_HighRank},
    $\Delta' := \col_{\cal C,\gamma}\Delta$ is a prime invariant subvariety,
    so it is a quasidiagonal in $\cal M_1'\times \cal M_2'$.
    It remains to show that $\cal M_i'$ is full rank. Let $M= (M_1,M_2)$
    and $\cal C_i$ be the cylinders of $\cal C$ on $M_i$.
    Let $M'=(M_1',M_2') = \col_{\cal C,\gamma}M$.
    Without loss of generality let $\gamma \subset M_1$.
    By Corollary \ref{collapse-full-rank}, $\cal M'_1 = 
    \col_{\cal C_1,\gamma} \cal M_1$ is full rank invariant subavariety
    with marked points.
    By Lemma \ref{collapse-one-dim},
    $\cal M_2'$ is codimension at most $1$.
    If $\cal M_2' \subset \cal M_2$, $\cal M_2'$ is full rank by Lemma
    \ref{lem:codimension_one_IS}. Otherwise, some saddle connection
    $\gamma' \subset M_2$ collapses, so
    $\cal M_2' = \col_{\cal C_2,\gamma'}\cal M_2$, so by 
    Corollary \ref{collapse-full-rank}, $\cal M_2'$ is full rank.
\end{proof}

\begin{defn}
Let $\cal M$ be an invariant subvariety in a 
potentially multi-component stratum.
Let $M\in \cal M$ be generic and $\cal C_1,\cal C_2$ be two distinct
$\cal M$-parallel classes of cylinders.
Let $\gamma_i$ be a cross curve on a cylinder in $\cal C_i$.
Assume furthermore that the components of $\col_{\cal C_1,\gamma_1,\cal C_2,\gamma_2}M$ have no 
translation surface automorphisms
other than the identity.
Then, $(\cal M, M, \cal C_1, \cal C_2, \gamma_1, \gamma_2)$ is called a
\define{good diamond}. Note that this definition is a special case of a skew diamond defined in \cite[Section 5]{ApisaWright_Geminal}.
\end{defn}

The following is similar to \cite[Lemma 3.31]{ApisaWright_Diamonds}.
\begin{lemma}\label{find-diamonds}
    Let $\cal M$ be a prime invariant subvariety of rank at least $2$. Then, there is a $G_\delta$ set $U\subset \cal M$ such that for $M\in U$ there is a good diamond that contains $M$. Furthermore, if $\cal M$ is rank at least $3$, there is a $G_\delta$ set of $M$ such $M$ contains three cylinders $\cal C_1,\cal C_2,\cal C_3$ such that there is a good diamond containing any two $\cal C_i$.  
\end{lemma}
\begin{proof}
    By Lemma \ref{generic}, a $G_\delta$ set $U'$ of surfaces in $\cal M$ are generic.
    Now we show that we can find $k$ cylinder classes on any surface $M$ in a dense open set $V$.
    By definition, $\cal M_i = \clos{\proj_i(\cal M)}$ has rank $k$.
    By \cite[Theorem 1.10]{Wright_Cylinder}, a dense subset of the set of horizontally
    periodic surfaces in $\cal M_i$ has at least $k$ disjoint cylinder
    equivalence classes.
    There is an open
    subset around each of these points, where these cylinders persist and
    remain disjoint, so a dense set of surfaces in $\cal M_i$ has at
    least $k$ disjoint cylinder equivalence classes.
    By Lemma \ref{m-parallel-quasidiagonal}, if $\proj_iM$ contains $k$ distinct $\cal M_i$-parallel classes, then $M$ contains $k$ distinct $\Delta$-parallel classes.
    Since $\proj_i$ is a
    submersion, there is a dense set of surfaces $\proj_i^{-1}(M) \subset \cal M$
    that have $k$ disjoint cylinder equivalence classes. Thus, the set $U''$ of generic surfaces in $\Delta$ that contain at least $k$ distinct $\Delta$-parallel classes is a $G_\delta$ set.
    
    For each $M\in U''$ there is an open set $V_M$ around $M$ where the $k$ cylinder classes persist. For any two of the cylinder classes $\cal C_1,\cal C_2$, the set of surfaces in $\col _{\cal C_1,\gamma_1,\cal C_2,\gamma_2} \cal M$ that have nontrivial translation surface isomorphisms is a set of isolated points. Thus, there is a countable union of two dimensional spaces of surfaces in $\Delta$ where $\cal C_1,\cal C_2$ do not belong to a good diamond. Thus, after removing these sets for all possible pairs of $\Delta$-parallel classes from $V_M$ we are left with an open dense set of $V_M$ that satisfy the lemma. Taking a union of all of these sets gives the $G_\delta$ set $U$ that is the conclusion of the lemma.
\end{proof}

\begin{defn} \label{def:diagonal-marked-points}
    Let $\cal M$ be an invariant subvariety with marked points. The 
    \define{diagonal}
    $\cal D \subset \cal M\times \cal M$ is the set of surfaces $(M,N)$ such 
    that there is a translation surface isomorphism $M\to N$ that takes 
    marked points to marked points. The \define{antidiagonal}
    $\overline {\cal D} \subset \cal M\times \cal M$ is the set of surfaces 
    $(M,N)$ such that the $N = -\Id(M)$ and $-\Id$ takes marked points to
    marked points.
\end{defn}

\begin{lemma} \label{diagonal_diamond}
Let $\Delta\subset \cal M_1\times \cal M_2$ be an (equal area) quasidiagonal.
If for any good diamond 
$(\Delta, M, \cal C_1, \cal C_2, \gamma_1, \gamma_2)$,
both $\col_{C_i,\gamma_i}\Delta$,
$i=1,2$ are diagonals (resp. antidiagonals) up to rescaling
(see Remark \ref{areas}),
then $\Delta$ is a diagonal (resp. antidiagonal).
\end{lemma}

\begin{proof}
We prove the statement for diagonals as the statement for antidiagonals is 
similar.
Let $(\Delta, M,\cal C_1,\cal C_2,\gamma_1,\gamma_2)$ be any good diamond,
and let $M = (M_1,M_2)$.
Assume both $\col_{\cal C_i,\gamma_i}\Delta$, $i=1,2$ are diagonals up to
rescaling. 
Then, there are constants $r_i>0$ and isomorphisms
$f_i:\col_{\cal C_i,\gamma_i} M_1 \to r_i\col _{\cal C_i,\gamma_i}M_2$ that restrict to isomorphisms
\begin{align*}
    &\col_{\cal C_2,\gamma_2}f_1: \col_{\cal C_1,\gamma_1,\cal C_2,\gamma_2}M_1
    \to r_1\col_{\cal C_1,\gamma_1,\cal C_2,\gamma_2}M_2\\
    &\col_{\cal C_1,\gamma_1}f_2: \col_{\cal C_1,\gamma_1,\cal C_2,\gamma_2}M_1
    \to r_2\col_{\cal C_1,\gamma_1,\cal C_2,\gamma_2}M_2
\end{align*}
And these must be the same isomorphism since we chose a good diamond. We also get that $r_1=r_2$ because isomorphic translation surfaces must have the
same area.
We can now define an isomorphism 
$f:M_1\to r_1M_2$ as follows.
$M_i-\overline{\cal C_j}$ can be identified with 
$\col_{\cal C_j,\gamma_j}M_i - \col_{\cal C_j,\gamma_j} \cal C_j$ for 
$i,j=1,2$.
Thus, $f$ can be defined to be $f_i$ on $M_1 - \overline{\cal C_i}$.
$f$ is defined twice on $M_1 - \overline{\cal C_1}-\overline{\cal C_2}$,
but these definitions agree.
This defines $f$ as an isomorphism of punctured translation surfaces
$M_1 -\overline{\cal C_1}\cap \overline{\cal C_2}
\to r_1M_2 -r_1\overline{\cal C_1}\cap \overline{\cal C_2}$.
$\overline{\cal C_1}\cap \overline{\cal C_2}$ consists of a finite set of 
points, so we can extend $f$ to an isomorphism from $M_1$ to $M_2$.

We note that $r_1=1$ because we assumed $M_1,M_2$ have the same area. 
By Lemma \ref{find-diamonds} for a dense set of $M$,
there exists a good diamond that contains $M$, and for these $M$ both
components $M_i$ are isomorphic,
so $\Delta$ is a diagonal.
\end{proof}

\subsection{Marked Points}

See \cite{ApisaWright_MarkedPoints} for a more in depth discussion of marked 
points.

\begin{defn} \label{def:marked-points}
    Let $\cal H$ be a stratum 
    and $\cal H^{*n}$ be the set of surfaces in
    $\cal H$ with $n$ distinct marked points. Let $\cal F:\cal H^{*n}\to \cal H$
    be the map that forgets marked points. Let $\cal M$ be an invariant 
    subvariety of $\cal H$. An \define{$n$-point marking} is an invariant subvariety
    of $\cal H^{*n}$ that maps to a dense subset of $\cal M$ under $\cal F$.
    We use the term point marking to refer to an $n$-point marking for any value of $n$.
    A point marking is \define{reducible} if it is a fiberwise union of two 
    other point markings, otherwise it is irreducible.
    A marked point $p$ on $M\in \cal M$ is $\cal M$-\define{free} if there are no relations in $\cal M$ between $p$ and any of the other marked points on $M$.
    When $\cal M$ is understood, we call $p$ a free marked point. Similarly when $\cal M$ is understood, we say a set of marked points on $M$ is irreducible if $\cal M$ is irreducible. 
\end{defn}

\begin{lemma} \label{lem:raise-genus}
    Let $\cal M$ be an invariant subvariety in a genus $g$ stratum. Let $M\in \cal M$ and $\cal C$ be an $\cal M$-parallel class that contains only a single simple cylinder $C$. Let $\gamma$ be a cross curve of $C$. Let $M' = \col_{\cal C,\gamma}M$ and $\cal M' =\col_{\cal C,\gamma}\cal M$. If the endpoints of $\col_{\cal C,\gamma}C$ are distinct from each other, then the genus $g(\cal M) = g(\cal M')+1$.
    
    Similarly, assume $\cal C$ contains two simple cylinders and assume the four endpoints of $\col_{\cal C,\gamma}\cal C$ contains at least three distinct points. Then $g(\cal M) = g(\cal M')+2$.  
\end{lemma}
\begin{proof}
    Let $\cal M'$ be contained in a stratum with marked points $\cal H = \cal H(\kappa)$, where $\kappa$ is a tuple containing the orders of the singularities of the surfaces in $\cal H$. Let $|\kappa|$ be the sum of the entries of $\kappa$ and $s$ be the number of elements in $\kappa$. Note that $2g-2 = |\kappa|$ for any stratum $\cal H(\kappa)$. The sum of the angles around all singularities of $M'$ is $2\pi(s+|\kappa|)$. After gluing in a cylinder, the sum of the angles around the singularities of $M$ is $2\pi(s+|\kappa|+1)$. But the two distinct endpoints of $\col_{\cal C,\gamma}C$  fuse into one, so $M$ has $s-1$ singularities, so the sum of the singularities is $|\kappa|+2 = 2g(\cal M) -2$. Thus, $g(\cal M) = g(\cal M')+1$. To prove the second statement, we can glue in the cylinders one at a time. Each time we increase the genus by one. 
\end{proof}

\begin{thm} \label{thm:marked-points}
    Let $\cal P$ be a nonempty irreducible point marking on 
    an invariant subvariety $\cal M$. If $\cal M$ is a full rank, then either $\cal P = \cal M^*$ or $\cal M$
    is a hyperelliptic locus with hyperelliptic involution $J$
    and $\cal P$ is one of the following
    \begin{align*}
        &\{(M, p): M\in \cal M,p \text{ is a Weierstrass point}\}\\
        &\{(M,p,J(p)): M\in \cal M,p\in M\}
    \end{align*}
    
\end{thm}

This theorem was proven for strata in \cite[Theorem 1.5]{Apisa_MarkedPoints}
and \cite[Theorem 1.4]{ApisaWright_MarkedPoints} in general.

\begin{lemma} 
    \label{lem:irreducible}
    Let $\cal M$ be an invariant subvariety with marked points. Let $M\in \cal M$,
    and let $\Gamma$ be a set of $\cal M$-parallel saddle connections on 
    $M$. Assume that no saddle connects a marked point to itself. Let $\gamma\in \Gamma$ be a saddle such that at least 
    one endpoint of $\gamma$ is a marked point but not a periodic point.
    Then, for each 
    endpoint $p$ of $\gamma$ that is a marked point but not a periodic point,
    there exists an irreducible set $P\ni p$ of marked points such that 
    for every $\gamma_i \in \Gamma$, $P$ contains an endpoint of $\gamma_i$.
\end{lemma}
\begin{proof}
    Let $P$ be a maximal 
    irreducible set of marked points of $M$ containing $p$. 
    We may move the points of $P$ without moving the rest of the surface.
    No saddle $\gamma'\in \Gamma$ can exist that does not have an endpoint in $P$
    otherwise we can make $\gamma$ not parallel to $\gamma'$, which is a 
    contradiction.
\end{proof}

\begin{lemma} \label{marked-points-from-collapse}
    Let $\cal M$ be an invariant subvariety (without marked points),
    $M\in \cal M$,
    $\cal C$ an $\cal M$-parallel class of simple
    cylinders on $M$, and $\gamma$ a cross curve of $\cal C$.
    Assume there is at least one marked point on 
    $M' := \col_{\cal C,\gamma}M$ and $\Gamma := \col_{\cal C,\gamma}\cal C$.
    \begin{enumerate}
        \item If $\cal M$ is a stratum, each marked point of $M'$ is free.
        Furthermore, $\Gamma$ contains a single saddle connection.

        \item If $\cal M$ is a hyperelliptic locus, and let $J$ be the induced hyperelliptic involution
        on $\cal M'$.
        Then, $M'$ contains no periodic points and $\Gamma$ either consists of one saddle connection fixed by $J$ or two saddle connections swapped by $J$. 
    \end{enumerate}
\end{lemma}
\begin{proof}
    First assume $\cal M$ is a stratum. By Theorem \ref{thm:boundary},
    every marked point is free. Let $\Gamma:=\col_{\cal C,\gamma}\cal C$. Let $\gamma\in \Gamma$ be a saddle such 
    that at least one endpoint of $\gamma$ is a marked point. 
    By Lemma \ref{lem:irreducible} and
    since all marked 
    points are free, $\gamma$ cannot be $\col_{\cal C,\gamma}\cal M$-parallel
    to any other saddle connection on $M'$. Thus, $\cal C$ only 
    consisted of a single cylinder, so there is only one saddle in $\Gamma$.
    Now assume $\cal M$ is a hyperelliptic locus.
    \begin{claim}
        None of the marked points on $M'$ are Weierstrass points.
    \end{claim}
    \begin{proof}
        Assume by contradiction that there exists a saddle $\gamma_1\in \Gamma$
    such that one of the endpoints is a Weierstrass point $p$.
    Let $\gamma_2 = J'(\gamma_1)$.
    We note that $\gamma_2 \neq \gamma_1$ otherwise $M$ has marked points,
    which contradicts our assumptions.
    Then, there are two cylinders 
    $C_1,C_2\in \cal C$, which  collapse to become $\gamma_1,\gamma_2$ respectively.
    Let $q_1$ be the other endpoint of $\gamma$, so $q_2 = J'(q_1)$ is the 
    other endpoint of $\gamma_2$. First assume $q_1 \neq q_2$.
    Note that $M$ is topologically $M'$ plus two handles, so it has genus 
    two more than $M'$.
    Thus, the hyperelliptic involution on $M$ must have four more fixed 
    points than the hyperelliptic involution on $M'$.
    By our assumptions, $J$ swaps $C_1$ and $C_2$.
    The fixed point of $J$ are the fixed points of $J'$
    minus $p$ plus the singularity that comes out of the fusion of $q_1,q_2$.
    Thus, $J$ same number of fixed points as $J'$.
    Thus, this is not possible.
    The remaining case is that $q_1=q_2$. In this case, $M$ has genus one 
    more than $M'$. However, there are two fewer fixed points of $J$
    than there are 
    fixed points of $J'$, which is also a contradiction. Thus, 
    $M'$ cannot have Weierstrass points.
    \end{proof}
    By the Claim and Theorem \ref{thm:marked-points}, $M'$ has no periodic points.
    $\cal C$ consists of at most two cylinders. If $\cal C$ is one cylinder, it must be fixed by the hyperelliptic involution. In this case, $\Gamma$ is one saddle $\gamma$, and $\gamma$ is fixed by $J$. Now assume that $\cal C$ contains two cylinders. These cylinders must be swapped by the hyperelliptic involution, so $\Gamma$ consists of two saddles swapped by $J$. This concludes the proof.
\end{proof}

\section{Genus 2}
\begin{figure}
    \centering
    \includegraphics[width=60mm]{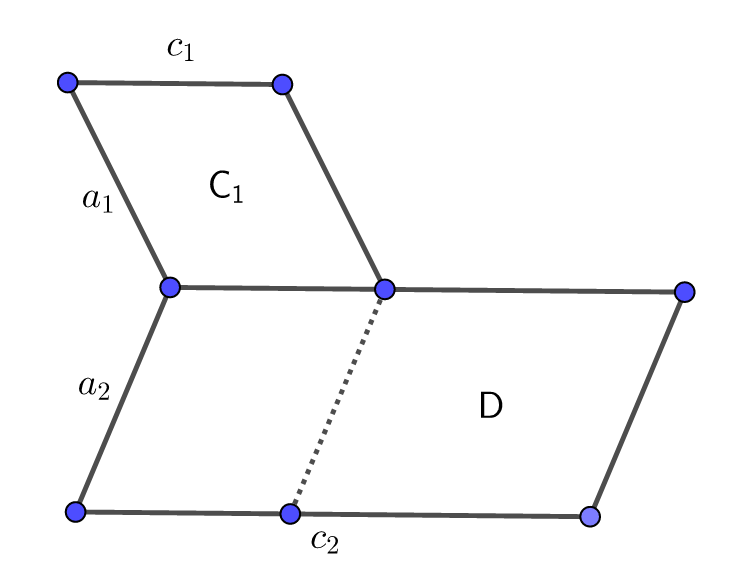}
    \captionof{figure}{Horizontally periodic surface is $\cal H(2)$.}
    \label{H(2)-figure}
\end{figure}

As a base case for the induction, we must prove Theorem \ref{main_theorem} for quasidiagonals in $\cal H(2)\times \cal H(2)$. 

\begin{thm}\label{theorem-H(2)}
The only (equal area) quasidiagonal $\Delta\in \cH(2)\times \cH(2)$ is the diagonal $\{(M,M):M\in \cH(2)\}$.
\end{thm}

Let $\Delta\subset \cH(2)\times \cH(2)$ be a quasidiagonal. Let $(M, M')\in \Delta$ be any generic surface. To prove the theorem, it suffices to show that $M'$ is equal to $M$, which is Lemma \ref{cor1} below.

\begin{lemma} \label{one-each-side}
    Let $\cal H_1,\cal H_2$ be hyperelliptic components, and $\Delta\subset \cal H_1,\times \cal H_2$ a quasidiagonal. Let $\cal C$ be a cylinder equivalence class on a surface $M\in \Delta$. Then, $\cal C$ consists of one cylinder on each component.
\end{lemma}
\begin{proof}
    By Lemma \ref{m-parallel-quasidiagonal}, $\cal C$ consists of one $\cal H_1$-parallel class and one $\cal H_2$-parallel class. For any stratum $\cal H$, a $\cal H$ parallel class consists of an equivalence class of homologous cylinders, and on a hyperelliptic component no two cylinders can be homologous.
\end{proof}

Choose two disjoint cylinders $C_1,D$ on $M$. By Lemma \ref{one-each-side}, there are cylinders $C_1',D'$ on $M'$ such that $\cal C_1= \{C_1,C_1'\}$ and $\cal D = \{D,D'\}$ are $\Delta$-parallel classes of cylinders. By Corollary \ref{cylinder_proportion}, $C_1',D'$ must be disjoint. 
We rotate $(M,M')$ to make $\cal C$ horizontal and perform a cylinder shear on $\cal D$ until $M$ is horizontally periodic. By Lemma \ref{horizontally_periodic}, $M'$ is also horizontally periodic. Let $C_2$ be the other horizontal cylinder on $M$. We fix a cross curve $\gamma_i$ of $C_i$, where $\gamma_2$ is a boundary curve of $D$. Let $a_i$ be the period of $\gamma_i$. Let $c_i$ be the period of the core curve of $C_i$. See Figure \ref{H(2)-figure}. Because $M$ is horizonally periodic, we have that $\im c_i=0$. Label the corresponding cylinders, saddle connections, and periods of $M'$ with primes.

Since $\Delta$ is cut out by linear equations in period coordinates, we have that
\[
T'\cdot \begin{pmatrix}
a_1'\\ a_2'\\ c_1'\\ c_2'
\end{pmatrix} = 
T\cdot 
\begin{pmatrix}
a_1 \\ a_2 \\ c_1 \\ c_2
\end{pmatrix}
\]
for some real matrices $T,T'$.
By Theorem \ref{ChenWright}, the absolute periods determine each other, so the above matrices are invertible. Thus we can assume $T'=\Id$. By Lemma \ref{drop_rank_lemma}, we may choose $a_1,a_1'$ so that they are $\Delta$-parallel. In addition, $a_2,a_2'$ are $\Delta$-parallel since they are boundary curves of $\Delta$-parallel cylinders. This means that $a_i'$ doesn't depend on any period except $a_i$. Changing $\im a_1,\im a_2$ does not affect $\im c_1',\im c_2'$ since the surface must remain horizontally periodic by Lemma \ref{horizontally_periodic}. Thus, we can simplify the matrix
\[
T 
=\begin{pmatrix}
f_{11} & 0 & 0 & 0\\
0 & f_{22} & 0 & 0\\
0 & 0 & g_{11} & g_{12}\\
0 & 0 & g_{21} & g_{22}
\end{pmatrix}.
\]

\begin{lemma} \label{g}
With the above notation, $g_{12}=g_{21}=0$.
\end{lemma}
\begin{proof}
We can dilate $\cal D$ while keeping the surface horizontally periodic. This changes the circumferences $c_2$ without changing $c_1$, so $g_{12} = 0$. 
The following combination of shears will change $c_1$ without changing $c_2'$ or $a_1'$.
Dilate $\cal D$ while keeping the surface horizontally periodic. Now shrink the whole surface in the real direction (i.e. by a matrix of the form
$\begin{pmatrix}
t & 0 \\
0 & 1
\end{pmatrix}$) so that $c_2$ is its original size. This changes $c_1$ without changing $c_2$, so $g_{21} = 0$. 
\end{proof}

\begin{lemma} \label{lem:same-moduli}
    Let $(M,M')$ be a generic surface in $\Delta$. 
    Let $C,C'$ be $\Delta$-parallel simple cylinders on $M,M'$ respectively. Then, they have the same modulus.
\end{lemma}
\begin{proof}
    Let $\cal C = \{C,C'\}$. Shear $\cal C$ so that a cross curve $\gamma$ of $C$ is vertical. By Lemma \ref{drop_rank_lemma}, $C'$ also has a vertical saddle connection $\gamma'$. Now continue to shear $\cal C$ until the first time $M$ once again has a vertical saddle connection. We notice that both must have been sheared one full rotation. Thus, the moduli must be equal.
\end{proof}

\begin{lemma}\label{cor1}
Let $\Delta\subset \cal H(2)\times\cal H(2)$ be a quasidiagonal and $(M,M')\in \Delta$ be a generic surface. Then $M,M'$ are isomorphic translation surfaces.
\end{lemma}

\begin{proof}
We label the surfaces using the notation above. Because the combinatorics of the surfaces are the same, it suffices to show that $T = \Id$. We first show $g_{22}=f_{22}$ and $f_{11}=g_{11}$ is similar. We take a sequence of generic surfaces in $\Delta$ that converge to $(M,M')$. For surfaces close enough to $(M,M')$ there are cylinders corresponding to $C_i,C_i'$. By Lemma \ref{lem:same-moduli},
$C_i,C_i'$ have the same modulus along this sequence, so they have the same moduli on $(M,M')$. 
Let $h_i,h_i'$ be the heights of $C_i,C_i'$ respectively. Then, $\frac{c_i}{h_i} = \frac{c_i'}{h_i'} = \frac{g_{ii}c_i}{f_{ii}h_i}$, so $f_{ii}=g_{ii}$.

Now we show $g_{11} = f_{22}$. We shear $\cal C_1,\cal C_2$ so that $\gamma_1,\gamma_2$ are vertical, so that $M,M'$ are both vertically and horizontally periodic. By the same argument as above, $D, D'$, must have the same modulus. Thus
\[
\frac{c_2-c_1}{a_2} = \frac{c_2'-c_1'}{a_2'} = \frac{g_{22}c_2-g_{11}c_1}{f_{22}a_2},
\]
so $f_{11}=f_{22}=g_{11}=g_{22}$. Now, we see that
\[
\begin{pmatrix}
a_1'\\ a_2'\\ c_1'\\ c_2'
\end{pmatrix} = 
f_{11}
\begin{pmatrix}
1 &&&\\
&1&&\\
&&1&\\
&&&1
\end{pmatrix}
\begin{pmatrix}
a_1 \\ a_2 \\ c_1 \\ c_2
\end{pmatrix}.
\]
Since we assumed both sides have the same area, we have that $f_{11}=1$. Thus, $T = \Id$. This finishes the proof of Theorem \ref{theorem-H(2)}.
\end{proof}

\section{Proof of Main Theorem}

By Theorem \ref{thm:full_rank}, Theorem \ref{main_theorem} is equivalent to the
following:

\begin{thm} \label{thm:main}
    Let $\cal M_1,\cal M_2$ be full rank invariant subvarieties (without marked points) in a stratum of translation surfaces in genus $g\ge 2$. There exists (equal area) quasidiagonals $\Delta \subset \cal M_1\times\cal M_2$ only if $\cal M_1 = \cal M_2$. In this case, $\Delta$ must be the diagonal and antidiagonal.
\end{thm}

\begin{proof}
    We use induction.
    Theorem \ref{theorem-H(2)} is the theorem for a quasidiagonal
    $\Delta\subset \cal H(2)\times \cal H(2)$, which is the base case.
    Let $r,s$ be given such that the theorem holds for any
    $\Delta \subset \cal M_1\times \cal M_2$ if either $\rk \Delta < r$ or $\rk \Delta = r$ and
    $\rel \cal M_1 +\rel \cal M_2<s$.
    Now, we will prove the theorem for $\rk \Delta = r$ and
    $\rel \cal M_1+\rel \cal M_2=s$.
    Let $M\in \Delta$ be a generic surface and $\cal C$ a $\Delta$-parallel class on $M$. Let $\gamma$ be a cross curve of a cylinder of $\cal C$.
    Let $\Delta' := \col_{\cal C,\gamma}\Delta $ and $\cal M_{j}' := \overline{\proj_j \Delta'}$.
    By Lemma \ref{full-rank-sides},
    $ \Delta' \subset \cal M_{1}'\times\cal M_{2}'$
    is a quasidiagonal, where $\cal M_{j}'$ are full rank invariant subvarieties
    with marked points.
    \begin{claim}
        $\Delta'\subset \cal M_{1}'\times \cal M_{2}'$
        is a diagonal or antidiagonal. Furthermore, if either $\cal M_1$ or $\cal M_2$ is a hyperelliptic locus, then $\Delta'$ is both a diagonal and antidiagonal.
        Recall that the diagonal and antidiagonal were defined for invariant
        subvarieties with marked points in Definition 
        \ref{def:diagonal-marked-points}.
    \end{claim}
    \begin{proof}
        When $\cal M_1'$ and $\cal M_2'$ both do not have marked points,
        then the claim is true by the induction hypothesis.
        Let $\cal F$ be the functor that forgets marked points. When marked 
        points are present, the induction hypothesis can still be used on
        $\cal F \Delta' \subset \cal F \cal M_1'\times \cal F\cal M_2'$.
        Define $\cal N := \cal F\cal M_1' = \cal F\cal M_2'$,
        which is full rank.
        Let $(M_1',M_2') = \col_{\cal C,\gamma}(M_1,M_2)$.
        Since $\cal F\Delta'$ is a diagonal or 
        quasidiagonal, there is a map between the surfaces 
        $I:\cal FM_1'\to \cal FM_2'$, which is either the identity map if $\cal F\Delta'$ is the diagonal or $-\Id$ if $\cal F\Delta'$ is the antidiagonal. When $\cal N$ is hyperelliptic, $\Delta'$ is a diagonal and an antidiagonal, so we let $I$ denote the identity and $I'$ denote $-\Id$.
        We may construct a translation surface with marked points $\widehat M$ by 
        taking the underlying surface $\cal FM_2'$, the marked points of 
        $M_2'$, and the image of the marked points of of $M_1'$ under $I$.
        By taking the set of all $\widehat M$, we get $\widehat{\cal N}$ an invariant subvariety with marked points.
        Let $\Gamma=\col_{\cal C,\gamma} \cal C$, defined in the beginning of Section \ref{sec:collapse}. As an abuse of notation, we will use to term saddle to include geodesic segments between marked points or singularities allowing marked points but not singularities on the interior of the segment. Let $\Gamma_1 = I(\Gamma \cap M_1')$ and $\Gamma_2 = \Gamma \cap M_2'$. Under our definition of saddle, $\Gamma_1$ and $\Gamma_2$ are sets of saddles on $\widehat{M}$.

        By Lemma \ref{marked-points-from-collapse}, none of the marked points are periodic points.
        In addition, no saddle of $\Gamma_1$ or $\Gamma_2$ can connect a marked point to itself because otherwise it would come from collapsing two marked points on $M$, but we assumed $M$ does not have marked points.
        Since $\Gamma$ is a $\Delta'$-parallel set of saddle connections on $M'$, then $\Gamma_1\cup \Gamma_2$ is a $\widehat {\cal N}$-parallel set of saddle connections on $\widehat{M}$. 
        Thus, by Lemma \ref{lem:irreducible}
        there is an irreducible set of marked points $P\subset \widehat{M}$ containing an endpoint of every saddle in $\Gamma_1\cup \Gamma_2$. By Theorem \ref{thm:marked-points}, $P$ consists of either one free marked point or two points that are swapped by the hyperelliptic involution. 
        
        \begin{subclaim} \label{sub}
            If $\gamma_1,\gamma_2\in \Gamma_1\cup \Gamma_2$ are $\widehat{\cal N}$-parallel saddles that share an endpoint that is also a marked point, then $\gamma_1 = \gamma_2$. 
        \end{subclaim}
        \begin{proof}
            Since $\gamma_1,\gamma_2$ don't contain any singularities (other than marked points) in the interior, there is a flat coordinate chart of $\widehat{M}$ that contains $\gamma_1$ and $\gamma_2$. Let $p$ be a marked point on $\gamma_1$ and $\gamma_2$.
            By Lemma \ref{marked-points-from-collapse}, no marked points of $\widehat{M}$ are Weierstrass point.
            Let $P$ be the maximal irreducible set of marked points containing $p$. 
            $\cal F\widehat{N}$ is a stratum or a hyperelliptic locus, so by Theorem \ref{thm:marked-points}, $P$ is either $\{p\}$ or $p,J(p)$, where $J$ is the hyperelliptic involution on $\widehat M$.
            We may perturb $P$ without changing the rest of the surface.
            If $P = \{p\}$, it is clear that if $\gamma_1\neq \gamma_2$, there is some perturbation that makes them not parallel. Since $\gamma_1,\gamma_2$ are $\widehat{\cal N}$-parallel, we must have $\gamma_1 = \gamma_2$. 
            Now assume $P = \{p,J(p)\}$. If $J(p)$ is not among the four endpoints of $\gamma_1,\gamma_2$, this is equivalent to the previous case. 
            We say that $\gamma_1,\gamma_2$ are on opposite sides of $p$ if the vectors from $p$ pointing towards the interior of $\gamma_i$ are 180 degrees from each other.
            If $\gamma_1$ and $\gamma_2$ are on opposite sides of $p$ then a small perturbation of $p$ will make them not parallel.
            If $\gamma_1,\gamma_2$ both connect $p$ and $J(p)$ and $\gamma_1\neq \gamma_2$, then $\gamma_1$ and $\gamma_2$ are on opposite sides of $p$, so perturbing will also cause $\gamma_1$ and $\gamma_2$ to not be parallel. The last case is that $\gamma_1$ connects $p$ and $J(p)$, and $\gamma_2$ connects $p$ to another point $q$. We perturb $p$ be a small vector $\delta$, so $J(p)$ is changed by $-\delta$. Thus, we change the period of $\gamma_1$ by $2\delta$ and the period of $\gamma_2$ by $\delta$. $\gamma_1$ and $\gamma_2$ are still parallel under any such perturbation, so the period of $\gamma_2$ must be half the period of $\gamma_1$. In addition, $\gamma_1$ and $\gamma_2$ are on the same side of $p$. Thus, $q$ is the midpoint of $\gamma_1$. However, this is a Weierstrass point since $\gamma_1$ is fixed by $J$, and no marked points of $\widehat{M}$ are Weierstrass points. This rules out the final case, so the subclaim is proven.
        \end{proof}

        By Subclaim \ref{sub}, each point in $P$ is adjacent to exactly one saddle. 
        Thus, $\Gamma_1\cup \Gamma_2$ consists of either one saddle connection or, in the case that $\cal N$ is a hyperelliptic locus two saddles. In the latter case by applying Subclaim \ref{sub} again, we get that the two saddles must be swapped by the hyperelliptic involution $J$.
        Now, using this understanding of $\Gamma_1\cup \Gamma_2$ and Lemma \ref{marked-points-from-collapse}, we can understand $\Gamma_1$ and $\Gamma_2$.

        \noindent \underline{Case 1:} $\cal M_1$ and $\cal M_2$ are strata.
        
        By Lemma \ref{marked-points-from-collapse}, $\Gamma_1,\Gamma_2$
        each contains a 
        single saddle $\gamma_1,\gamma_2$ respectively.
        Thus, $\gamma_1$ and $\gamma_2$ either must be the same saddle connection, or they are swapped by the hyperelliptic involution.
        If $\gamma_1=\gamma_2$, $I$ takes 
        marked points to marked points, so $\Delta'$ is a diagonal or 
        antidiagonal. Otherwise $\cal N$ is hyperelliptic and $I'$ takes marked points to marked points, so $\Delta'$ is an antidiagonal.

        \noindent\underline{Case 2:} 
        $\cal M_1$ and $\cal M_2$ are a hyperelliptic loci.

        If $\Gamma_1\cup \Gamma_2$ consists of a single saddle connection, then $\Gamma_1=\Gamma_2$ is a single saddle and by Lemma \ref{marked-points-from-collapse} this saddle is fixed by the hyperelliptic involution. Thus both $I,I'$ take marked points to marked points, so $\Delta'$ is both a diagonal and antidiagonal.
        Now assume $\Gamma_1 \cup \Gamma_2 = \{\gamma,J(\gamma)\}$, where $\gamma \neq J(\gamma)$.
        By Lemma \ref{marked-points-from-collapse},
        $\Gamma_i$ contains one saddle that is fixed by $J$ or contains 
        two saddles that are swapped by $J$. Since none of the saddles in $\Gamma_1\cup \Gamma_2$ are fixed by $J$, $\Gamma_1=\Gamma_2 = \{\gamma,J(\gamma)\}$, and $\Delta'$ is both a diagonal and antidiagonal.

        \noindent \underline{Case 3:}
        $\cal M_1$ is a stratum and $\cal M_2$ is a hyperelliptic locus.

        By Lemma \ref{marked-points-from-collapse}, $\Gamma_1$ is a single saddle connection, and $\Gamma_2$ can be either one saddle fixed by $J$ or two saddle connections swapped by $J$.
        We will show be contradiction that $\Gamma_2$ cannot consists of two distinct saddle connections.
        Assume $\Gamma_2 = \{\gamma,J(\gamma)\}$, $\gamma \neq J(\gamma)$.
        First we will show that among the four endpoints of $\gamma,J(\gamma)$, there are at least three distinct points.
        Using Lemma \ref{lem:irreducible} on $\widehat{M}$, $\gamma$ and $J(\gamma)$ must each have at least one marked point as an endpoint. 
        For each saddle, its two endpoints cannot be the same marked point because we assumed $\cal M_1$ and $\cal M_2$ do not have marked points. In addition, the marked points of $\gamma$ and $J(\gamma)$ are distinct by the Subclaim. Thus, among the four endpoints there are at least two distinct marked points along with at least one more distinct endpoint. Thus by Lemma \ref{lem:raise-genus}, the genus $g(\cal M_1) \le g(\cal N) +1$ but $g(\cal M_2) = g(\cal N)+2$, which is a contradiction since we assumed $\cal M_1$ and $\cal M_2$ are full rank invariant subvarieties with the same rank. Thus, $\Gamma_2$ can only have a single saddle which is fixed by $J$.

        Recall $\Gamma_1\cup\Gamma_2$ is either a single saddle fixed by $J$ or two saddles swapped by $J$.
        But $\Gamma_1\cup \Gamma_2$ cannot be two saddles swapped by the hyperelliptic involutions because $\Gamma_2$ contains a saddle fixed by $J$. Thus, 
        $\Gamma_1 \cup \Gamma_2$ must consist of only a single saddle connection, so $\Gamma_1=\Gamma_2$, and we get that $\Delta'$ is both a diagonal and antidiagonal. The case where $\cal M_1$ is a hyperelliptic locus and $\cal M_2$ is a stratum is equivalent to this case, so we have proven the claim.
    \end{proof}
    
    Let $(\cal M, M, \cal C_1, \cal C_2, \gamma_1,\gamma_2)$ be a good diamond, which exists by Lemma \ref{find-diamonds}.
    By the claim, if $\cal M_1$ or $\cal M_2$ is hyperelliptic, then $\col _{\cal C_i, \gamma_i}\Delta$ is a diagonal for $i=1,2$.
    Thus by Lemma \ref{diagonal_diamond},
    $\Delta$ must be a diagonal. Now it remains to consider when
    $\cal M_1$ and $\cal M_2$ are not hyperelliptic.
    In this case, the rank is at least $3$,
    so by Lemma \ref{find-diamonds} we can find $M$ and three disjoint equivalence
    classes of cylinders $\cal C_1,\cal C_2,\cal C_3$ such that for any two $\cal C_i\neq \cal C_j$, we can find a good diamond $(\cal M, M, \cal C_i,\cal C_j,\gamma_i,\gamma_j$). By the Claim for $i=1,2,3$, 
    $\col _{\cal C_i,\gamma_i}\Delta$ is a diagonal or quasidiagonal.
    By the Pigeon Hole Principle,
    there are two $(\cal C_i,\gamma_i)$ such that
    $\col _{\cal C_i,\gamma_i}\Delta$ are either both diagonals or both
    quasidiagonals. Now Lemma \ref{diagonal_diamond} finishes the proof.
\end{proof}

\section{Appendix}
In this section, we prove Theorem \ref{SmillieWeiss}, 
which is a multicomponent version of \cite[Corollary 6]{SmillieWeiss}. 
Most of the work will be 
adapting \cite[Corollary 2.7]{MinskyWeiss} to the multicomponent setting.
Let $\cal Q$ be a stratum of
multicomponent quadratic differentials and $\cal Q_1$ the unit area locus of
$\cal Q$. Define $h_t = 
\begin{pmatrix}
    1 & t\\
    0 & 1
\end{pmatrix} \in \SL(2,\R)$, and $H_t = \{h_t,t\in \R\}$.

\begin{thm} \label{minimal-sets}
    Every closed $H_t$-invariant set of $\cal Q_1$
    contains a minimal closed $H_t$-invariant
    set. A minimal closed $H_t$-invariant set is compact.
\end{thm}

The proof of Theorem \ref{minimal-sets} will be broken down into Proposition
\ref{minimal-sets-exist} and Proposition \ref{minimal-sets-compact}, which
depend on the following nondivergence result. Let 
\[
    \Avg_{T,q}(K) := \frac{|\{t\in[0,T]:u_tq\in K\}|}{T}
\]

\begin{thm}\label{nondivergence}
    For any $\epsilon>0$ and $\eta>0$, there is a compact subset $K$
    such that for any
    $q\in \cal Q_1$, one of the following statements holds:
    
    \begin{enumerate}
        \item $\displaystyle \liminf_{T\to\infty} 
        \Avg_{T,q}(K) \ge 1-\epsilon$.
        \item $q$ contains a horizontal saddle connection of length less than
        $\eta$.
    \end{enumerate}
\end{thm}

Theorem \ref{nondivergence} is very similar to but more general than 
Theorem H2 from \cite{MinskyWeiss}. 
Theorem \ref{nondivergence} is not stated in this way
in \cite{MinskyWeiss}, but this formulation is more coninvient for our purposes.
The proof of this theorem follows 
exactly the proof of Theorem H2 in the single component case, and
we copy it here for the coninvience of the reader. It mainly uses Theorem 6.3
of \cite{MinskyWeiss}, which we prove for the multicomponent case,
assuming the single component case.

Let
$\widetilde{\cal Q}$ be the marked quadratic
differentials that cover $\cal Q$. Let $\widetilde{\cal Q}_1$
be the unit area locus of $\widetilde{\cal Q}$. 
For $q\in \widetilde{\cal Q}$, let
$\cal L_q$ be the set of saddle connections on $q$. Let $l_{q,\delta}(t)$ be
the length of the saddle $\delta$ on $u_tq$, and $\|\cdot \|_I$ the 
$L^\infty$-norm on an interval $I\subset \R$. Let $\alpha_q(t)$ be the length of the 
shortest saddle connection on $u_tq$.

\begin{thm} \label{nondivergence-helper}
    There are positive constants $C,\alpha,\rho_0$, depending only on the 
    genus,
    such that if $q\in
    \widetilde{\cal Q}_1$, an interval $I\subset \R$, and $0<\rho'\le \rho_0$,
    satisfy:
    \[
        \text{for any } \delta\in \cal L_q, \|l_{q,\delta}\|_I \ge \rho',
    \] 
    then for any $0<\epsilon<\rho'$ we have:
    \[    
        |\{t\in I :\alpha_q(t)<\epsilon\}| 
        \le C \left(\frac{\epsilon}{\rho'} \right)^\alpha |I|
    \]
\end{thm}
\begin{proof}
    Let $q = (q_1,\dots,q_n)$ be a multicomponent quadratic differential. The
    condition $\forall\delta\in \cal L_q, \|l_{q,\delta}\|_I \ge \rho'$ implies
    this condition holds on each component. By the single component version of
    this theorem 
    \cite[Theorem 6.3]{MinskyWeiss}, $|\{t\in I :\alpha_{q_i}(t)<\epsilon\}| 
    \le C_i \left(\frac{\epsilon}{\rho'} \right)^{\alpha_i} |I|$
    and for some
    $\rho' < \min(\rho_{0,1},\dots,\rho_{0,n})$ for all $i$. Thus,
    \[
        |\{t\in I :\alpha_{q}(t)<\epsilon\}|
        = 
        |\bigcup_i \{t\in I :\alpha_{q_i}(t)<\epsilon\}|
        \le n\max_i C_i
        \left(\frac{\epsilon}{\rho'} \right)^{\min_i\alpha_i} |I|
    \]
\end{proof}

\begin{proof}[Proof of Theorem \ref{nondivergence}]
    Let $C,\alpha,\rho_0$ be the constants from Theorem 
    \ref{nondivergence-helper}.
    For fixed $\epsilon$ and $\eta$, choose $\epsilon'$ small enough so that 
    \[
        C\left(\frac{\epsilon'}{\rho_0}\right)^\alpha< \epsilon \text{, }
        \epsilon' < \rho_0 \text{, and } \epsilon' < \eta.
    \]
    Let $K$ be the sets of surfaces in $\cal Q_1$ with no saddle of length less
    than $\epsilon'$. This is compact by Masur's Compactness Criterion. 
    Let $q\in \widetilde{\cal Q}$, and suppose $q$ does not contain a
    horizontal saddle of length less than $\eta$. 
    Let $\rho = \min(\rho_0,\eta)$. The set
    \[
        \cal L_0 = \{\delta \in \cal L _q: l_{q,\delta}(0)<\rho\}
    \]
    is finite by \cite[Proposition 4.8]{MinskyWeiss}. Since we are assuming 
    none of the functions $t\mapsto l_{q,\delta}(t)$ are constant for $\delta
    \in \cal L_0$, they diverge by \cite[Lemma 4.4]{MinskyWeiss}. Thus there
    is some $T_0$ such that for all 
    $\delta\in \cal L_0, l_\delta(T_0) \ge \rho$. 
    For any $T\ge T_0$ we can apply Theorem \ref{nondivergence-helper} with
    $I = [0,T]$ and $\rho' = \rho$, and obtain that
    \[
        \Avg_{T,q}(\cal Q_1- K ) <\epsilon.
    \]
\end{proof}

\begin{prop} \label{minimal-sets-exist}
    Let $\cal Q_1$ be a stratum of unit-area
    multi-component quadratic differentials.
    Every closed $H_t$-invariant set $X\subset \cal Q_1$ contains a minimal 
    closed $H_t$-invariant set.     
\end{prop}

\begin{lemma} \label{const-functional-lemma}
    Let $X$ be a closed invariant set such that
    \[
        \rho := \inf\{l_{q,\delta}(0): q\in X,
        \delta \text{ is a horizontal saddle connection on }q \} > 0.
    \]
    Then $X$ contains a minimal closed invariant set. 
\end{lemma}

\begin{proof}
    Choose any $0<\epsilon<1$ and $0 < \eta < \rho$ 
    and let $K$ be the compact set obtained from Theorem \ref{nondivergence}.
    Then $K$ intersects $U_tq$ for all
    $q\in X$. Let $\{X_\alpha\}$ be any totally ordered family of 
    closed invariant subsets of $X$. Any finite intersection is nonempty
    $K\cap X_{\alpha_1}\cap \cdots \cap X_{\alpha_i}$ since every $H_t$ 
    trajectory meets $K$. Thus, $K \cap \bigcap_\alpha X_\alpha$ is nonempty
    by compactness. 
    By Zorn's lemma, $X$ contains a minimal closed invariant set.
\end{proof}

\begin{proof}[Proof of Proposition \ref{minimal-sets-exist}]
    We must prove that $X$ contains a set $X_0$ that satisfies the hypotheses
    of Lemma \ref{const-functional-lemma}. Let $q\in X$ have the maximal
    number of horizontal saddle connections, and let $X_0 = \clos{U_tq}$.
    For any $q'\in X_0$, the set of horizontal saddles on $q$ is isometric
    to a set of horizontal saddles on $q'$.
    Since $q$ has the maximal number
    of horizontal saddles, it cannot have more horizontal saddles. Thus,
    every surface in $X_0$ has the same horizontal saddles, so they have the
    same value of $\rho$.
\end{proof}

We do not include a proof of the following lemma 
which is Lemma 7.2 in \cite{MinskyWeiss}.

\begin{lemma} \label{general-compactness-lemma}
    Let $\{T_t\}$ be an action of $\bb R$ by homeomorphisms on a locally
    compact space $Z$. Suppose there is a compact $K\subset Z$
    such that for every $z\in Z$,
    the subsets $\{t\ge 0: T_tz\in K\}$ and $\{t\le 0: T_tz\in K\}$ are
    unbounded. Then $Z$ is compact. 
\end{lemma}

\begin{prop} \label{minimal-sets-compact}
    If a minimal $H_t$-invariant set exists, then it is compact.
\end{prop}

\begin{proof}
    By Lemma \ref{general-compactness-lemma} is suffices to prove that 
    $\{t\ge 0: u_tq\in X\}$ and $\{t\le 0: u_tq\in X\}$ are dense in $X$. We
    prove this for $\{t\ge 0: u_tq\in X\}$ as $\{t\le 0: u_tq\in X\}$ is 
    similar.
    Define $X^+(q)$ to be the limit points of $\{t\ge 0: u_tq\in X\}$, so it
    is a closed $H_t$-invariant set. It is nonempty by Theorem 
    \ref{nondivergence}. Since $X$ is minimal, $X^+(q) = X$. Thus, 
    $\{t\ge 0: u_tq\in X\}$ is dense in $X$.
\end{proof}

\begin{proof}[Proof of Theorem \ref{SmillieWeiss}]
    By Theorem \ref{minimal-sets}, a $H_t$-orbit closure contains a 
    minimal closed $H_t$-invariant set $X$, and $X$ is compact.
    Let $q = (q_1,\dots,q_n)\in X$, 
    so $\overline{U_tq} \subset X$. Let $\proj_i$ be the projection
    onto the $i$-th component. 
    $\clos{U_tq_i}$ is contained in a compact set, namely
    $\proj_i(X)$, so by \cite[Theorem 5]{SmillieWeiss}, 
    $q_i$ is horizontally periodic. This argument does not depend on $i$, so
    every component of $q$ is horizontally periodic.
\end{proof}

\bibliographystyle{halpha}
\bibliography{bib}

\begin{thebibliography}{EMM15}

\bibitem[AEM17]{AEM}
Artur Avila, Alex Eskin, and Martin M\"{o}ller.
\newblock Symplectic and isometric {${\rm SL}(2,\mathbb R)$}-invariant
  subbundles of the {H}odge bundle.
\newblock {\em J. Reine Angew. Math.}, 732:1--20, 2017.

\bibitem[Api20]{Apisa_MarkedPoints}
Paul Apisa.
\newblock {${\rm GL}_2\mathbb{R}$}-invariant measures in marked strata: generic
  marked points, {E}arle-{K}ra for strata and illumination.
\newblock {\em Geom. Topol.}, 24(1):373--408, 2020.

\bibitem[AW21a]{ApisaWright_HighRank}
Paul Apisa and Alex Wright.
\newblock High rank invariant subvarieties, 2021, arXiv:2102.06567.

\bibitem[AW21b]{ApisaWright_MarkedPoints}
Paul Apisa and Alex Wright.
\newblock Marked points on translation surfaces.
\newblock {\em Geom. Topol.}, 25(6):2913--2961, 2021.

\bibitem[AW21c]{ApisaWright_Diamonds}
Paul Apisa and Alex Wright.
\newblock Reconstructing orbit closures from their boundaries.
\newblock {\em Mem. Amer. Math. Soc.}, to appear 2021, arXiv:2011.08807.

\bibitem[AW22]{ApisaWright_Geminal}
Paul Apisa and Alex Wright.
\newblock Generalizations of the {E}ierlegende-{W}ollmilchsau.
\newblock {\em Camb. J. Math.}, 10(4):859--933, 2022.

\bibitem[CW21]{ChenWright_WYSIWYG}
Dawei Chen and Alex Wright.
\newblock The {WYSIWYG} compactification.
\newblock {\em J. Lond. Math. Soc. (2)}, 103(2):490--515, 2021.

\bibitem[EM18]{EM}
Alex Eskin and Maryam Mirzakhani.
\newblock Invariant and stationary measures for the {${\rm SL}(2,\mathbb R)$}
  action on moduli space.
\newblock {\em Publ. Math. Inst. Hautes \'{E}tudes Sci.}, 127:95--324, 2018.

\bibitem[EMM15]{EMM}
Alex Eskin, Maryam Mirzakhani, and Amir Mohammadi.
\newblock Isolation, equidistribution, and orbit closures for the {${\rm
  SL}(2,\mathbb R)$} action on moduli space.
\newblock {\em Ann. of Math. (2)}, 182(2):673--721, 2015.

\bibitem[Fil16]{Filip}
Simion Filip.
\newblock Splitting mixed {H}odge structures over affine invariant manifolds.
\newblock {\em Ann. of Math. (2)}, 183(2):681--713, 2016.

\bibitem[MW02]{MinskyWeiss}
Yair Minsky and Barak Weiss.
\newblock Nondivergence of horocyclic flows on moduli space.
\newblock {\em J. Reine Angew. Math.}, 552:131--177, 2002.

\bibitem[MW17]{MirzakhaniWright_Boundary}
Maryam Mirzakhani and Alex Wright.
\newblock The boundary of an affine invariant submanifold.
\newblock {\em Invent. Math.}, 209(3):927--984, 2017.

\bibitem[MW18]{MirzakhaniWright_Full}
Maryam Mirzakhani and Alex Wright.
\newblock Full-rank affine invariant submanifolds.
\newblock {\em Duke Math. J.}, 167(1):1--40, 2018.

\bibitem[MZ08]{MasurZorich}
Howard Masur and Anton Zorich.
\newblock Multiple saddle connections on flat surfaces and the principal
  boundary of the moduli spaces of quadratic differentials.
\newblock {\em Geom. Funct. Anal.}, 18(3):919--987, 2008.

\bibitem[NW14]{NW}
Duc-Manh Nguyen and Alex Wright.
\newblock Non-{V}eech surfaces in {$\cal{H}^{\rm hyp}(4)$} are generic.
\newblock {\em Geom. Funct. Anal.}, 24(4):1316--1335, 2014.

\bibitem[SW04]{SmillieWeiss}
John Smillie and Barak Weiss.
\newblock Minimal sets for flows on moduli space.
\newblock {\em Israel J. Math.}, 142:249--260, 2004.

\bibitem[Wri14]{Wright_Field}
Alex Wright.
\newblock The field of definition of affine invariant submanifolds of the
  moduli space of abelian differentials.
\newblock {\em Geom. Topol.}, 18(3):1323--1341, 2014.

\bibitem[Wri15]{Wright_Cylinder}
Alex Wright.
\newblock Cylinder deformations in orbit closures of translation surfaces.
\newblock {\em Geom. Topol.}, 19(1):413--438, 2015.

\end{thebibliography}

\end{document}